\documentclass[12pt]{article}

\usepackage{amsfonts}
\usepackage{latexsym}

\def\medskipamount{8pt} 
\def\arraystretch{1.3}
\begin{document} 


\renewcommand{\thesection}%
{Part \Roman{section}\hspace{.2ex}:\hspace{-1ex}}
\renewcommand{\thesubsection}{\arabic{subb}}

\newcounter{subb}
\setcounter{subb}{0}

\newcommand{\bit}[1]{\pagebreak[3] \addtocounter{subb}{1}
  \subsection{\hspace{-2ex} #1}\setcounter{equation}{0}\penalty12000}
\renewcommand{\theequation}{\thesubsection.\arabic{equation}}


\newcommand{\re}[1]{\mbox{\bf (\ref{#1})}}


\catcode`\@=\active
\catcode`\@=11
\def\@eqnnum{\hbox to .01pt{}\rlap{\bf \hskip -\displaywidth(\theequation)}}
\catcode`\@=12



\catcode`\@=\active
\catcode`\@=11
\newcommand{\nc}{\newcommand}

   
\nc{\bs}[1]{\addpenalty{-1000}
\addvspace{\medskipamount}\penalty10000 
\refstepcounter{equation}
\noindent \begin{em}{\bf (\theequation)} \nopagebreak}

\nc{\es}{\par \end{em} \addvspace{\medskipamount} } 

\nc{\br}[1]{\addvspace{\medskipamount} 
\refstepcounter{equation} 
\noindent {\bf (\theequation) #1.} \nopagebreak }

\nc{\er}{\par \addvspace{\medskipamount} }

\newcounter{index}
\nc{\bl}{\begin{list}{{\rm (\roman{index})}}{\usecounter{index}}}
   
\nc{\el}{\end{list}}

\nc{\pf}{\addvspace{\medskipamount} \par \noindent {\em Proof}}

\nc{\fp}{\phantom{.} \hfill \mbox{     $\Box$} \par
\addvspace{\medskipamount}} 

\nc{\beq}{\begin{equation}}
\nc{\eeq}{\end{equation}}

\nc{\beqas}{\begin{eqnarray*}}
\nc{\eeqas}{\end{eqnarray*}}


\nc{\C}{\mathbb C}
\nc{\Hyp}{\mathbf H}
\nc{\Pj}{\mathbb P}
\nc{\Q}{\mathbb Q}
\nc{\R}{\mathbf R}
\nc{\Z}{\mathbb Z}
 

\nc{\oper}[1]{\mathop{\mathchoice{\mbox{\rm #1}}{\mbox{\rm #1}}
{\mbox{\rm \scriptsize #1}}{\mbox{\rm \tiny #1}}}\nolimits}
\nc{\Ad}{\oper{Ad}}
\nc{\const}{\oper{const.}}
\nc{\ch}{\oper{ch}}
\nc{\End}{\oper{End}}
\nc{\END}{\oper{$\mathbf{End}$}}
\nc{\Hom}{\oper{Hom}}
\nc{\HOM}{\oper{$\mathbf{Hom}$}}
\nc{\id}{\oper{id}}
\nc{\Jac}{\oper{Jac}}
\nc{\Map}{\oper{Map}}
\nc{\mod}{\oper{mod}}
\nc{\Pol}{\oper{Pol}}
\nc{\Quot}{\oper{Quot}}
\nc{\slope}{\oper{slope}}
\nc{\Sym}{\oper{Sym}}
\nc{\rk}{\oper{rk}}
\nc{\td}{\oper{td}}
\nc{\tr}{\oper{tr}}

\nc{\operlimits}[1]{\mathop{\mathchoice{\mbox{\rm #1}}{\mbox{\rm #1}}
{\mbox{\rm \scriptsize #1}}{\mbox{\rm \tiny #1}}}}

\nc{\medoplus}{\operlimits{$\bigoplus$}}


\nc{\al}{\alpha}
\nc{\be}{\beta}
\nc{\ep}{\varepsilon}
\nc{\ga}{\gamma}
\nc{\Ga}{\Gamma}
\nc{\la}{\lambda}
\nc{\La}{\Lambda}
\nc{\si}{\sigma}
\nc{\Sig}{{\Sigma}}
\nc{\Om}{\Omega}


\nc{\elll}{l}
\nc{\tag}{}

\nc{\down}[1]{{\phantom{\scriptstyle #1}
\hbox{$\left\downarrow\vbox to
    9.5pt{}\right.\nulldelimiterspace=0pt \mathsurround=0pt$}
\raisebox{.4ex}{$\scriptstyle #1$}}}
\nc{\leftdown}[1]{
\hbox{$\raisebox{.4ex}{$\scriptstyle #1$}\left\downarrow\vbox to
    9.5pt{}\right.\nulldelimiterspace=0pt \mathsurround=0pt$}
    \phantom{\scriptstyle #1}}
\nc{\bigdowneq}{{\Big| \! \Big|}}

\nc{\ca}{{\mathcal A}}
\nc{\ce}{{\mathcal E}}
\nc{\ci}{{I}}
\nc{\cl}{{\mathcal L}}
\nc{\co}{{\mathcal O}}
\nc{\cu}{{\mathcal U}}
\nc{\cv}{{\mathcal V}}
\nc{\B}{{\mathcal B}}
\nc{\F}{{\mathcal F}}
\nc{\G}{{\mathcal G}}
\renewcommand{\H}{{\mathcal H}}
\nc{\M}{{\mathcal M}}
\nc{\N}{{\mathcal N}}

\nc{\GL}[1]{{{\rm GL(}#1{\rm )}}}
\nc{\PGL}[1]{{{\rm PGL(}#1{\rm )}}}
\nc{\SL}[1]{{{\rm SL(}#1{\rm )}}}
\nc{\SU}[1]{{{\rm SU(}#1{\rm )}}}
\nc{\U}[1]{{{\rm U(}#1{\rm )}}}
\nc{\Sp}[1]{{{\rm Sp(}#1{\rm )}}}
\nc{\cx}{\C^\times}

\nc{\bino}[2]{\mbox{\Large $#1 \choose #2$}}

\nc{\half}{\mathchoice{{\textstyle \frac{\scriptstyle 1}{\scriptstyle
2}}} {{\textstyle\frac{\scriptstyle 1}{\scriptstyle 2}}}
{\frac{\scriptscriptstyle 1}{\scriptscriptstyle 2}}
{\frac{\scriptscriptstyle 1}{\scriptscriptstyle 2}}}
\nc{\quarter}{\mathchoice{{\textstyle \frac{\scriptstyle
1}{\scriptstyle 4}}} {{\textstyle\frac{\scriptstyle 1}{\scriptstyle
4}}} {\frac{\scriptscriptstyle 1}{\scriptscriptstyle 4}}
{\frac{\scriptscriptstyle 1}{\scriptscriptstyle 4}}}
\nc{\ratio}[2]{\mathchoice{ {\textstyle \frac{\scriptstyle
#1}{\scriptstyle #2}}} {{\textstyle\frac{\scriptstyle #1}{\scriptstyle
#2}}} {\frac{\scriptscriptstyle #1}{\scriptscriptstyle #2}}
{\frac{\scriptscriptstyle #1}{\scriptscriptstyle #2}}}

\nc{\stack}[2]{
{\def\arraystretch{.7}\def\arraycolsep{0pt}\begin{array}{c}
\scriptstyle #1 \\ \scriptstyle #2 \end{array}\def\arraystretch{1}} }


\nc{\Left}[1]{\hbox{$\left#1\vbox to
    11.5pt{}\right.\nulldelimiterspace=0pt \mathsurround=0pt$}}
\nc{\Right}[1]{\hbox{$\left.\vbox to
    11.5pt{}\right#1\nulldelimiterspace=0pt \mathsurround=0pt$}}
\nc{\LEFT}[1]{\hbox{$\left#1\vbox to
    15.5pt{}\right.\nulldelimiterspace=0pt \mathsurround=0pt$}}
\nc{\RIGHT}[1]{\hbox{$\left.\vbox to
    15.5pt{}\right#1\nulldelimiterspace=0pt \mathsurround=0pt$}}

\nc{\lp}{\raisebox{-.1ex}{\rm\large(}\hspace{-.2ex}}
\nc{\rp}{\hspace{-.2ex}\raisebox{-.1ex}{\rm\large)}}

\nc{\emb}{\hookrightarrow}
\nc{\jay}{j}
\nc{\lrow}{\longrightarrow}
\nc{\llrow}{\hbox to 25pt{\rightarrowfill}}
\nc{\lllrow}{\hbox to 29pt{\rightarrowfill}}
\nc{\mbar}{\overline{M}}
\nc{\dbar}{\overline{\partial}}
\nc{\sans}{\,\backslash\,}
\nc{\st}{\, | \,}

\nc{\bee}{{\bf E}}
\nc{\bphi}{{\bf \Phi}}

\nc{\rs}{C}

\nc{\T}{T}
\nc{\Tor}{{\mathcal T}}

\nc{\g}{{\mathfrak g}}

\renewcommand{\labelenumi}{{\rm (\alph{enumi})}}

\catcode`\@=12



\noindent
{\LARGE \bf Generators for the cohomology ring 
\smallskip \\of the moduli space
of rank 2 Higgs bundles}\bigskip\\
{\bf Tam\'as Hausel }\\
Department of Mathematics, University of California, 
Berkeley, Calif. 94720\smallskip \\
{\bf Michael Thaddeus } \\
Department of Mathematics, Columbia University, 
New York, N.Y. 10027  
\renewcommand{\thefootnote}{}
\footnotetext{T.H. supported by NSF grant DMS--97--29992;
  M.T. supported by NSF grant DMS--98--08529.} 

\bigskip

\noindent
A central object of study in gauge theory is the moduli space of
unitary flat connections on a compact surface.  Thanks to the efforts of
many people, a great deal is understood about the ring structure of its
cohomology.  In particular, the ring is known to be generated by the
so-called {\em universal classes} \cite{ab,new}, and, in rank 2, all the
relations between these classes are also known \cite{bar,kn,st,zag}.

If instead of just unitary connections one allows {\em all\/} flat
connections, one obtains larger moduli spaces of equal importance and
interest.  However, these spaces are not compact and so very little
was known about the ring structure of their cohomology.

This paper will show that, in the rank 2 case, the cohomology ring of
this noncompact space is again generated by universal classes.  A
companion paper \cite{ht} gives a complete set of explicit relations
between these generators.

The noncompact spaces studied here have significance extending well
beyond gauge theory.  They play an important role in 3-manifold
topology: see for example Culler-Shalen \cite{cs}.  And they are the
setting for much of the geometric Langlands program: see for example
Beilinson-Drinfeld \cite{bd}.  But they have received perhaps the most
attention from algebraic geometers, in the guise of moduli spaces of
{\em Higgs bundles}.  A Higgs bundle is a holomorphic object, related to
a flat connection by a correspondence theorem similar to that of
Narasimhan-Seshadri in the unitary case.  This point of view has been
exploited to great effect, notably by Hitchin \cite{h,h2,h3} and Simpson
\cite{simp3,ubi,simp,simp2}.

The Higgs point of view predominates in this paper also.  Indeed, it
is strongly influenced by, and occasionally parallel to, the works of
Hitchin \cite{h} and Atiyah-Bott \cite{ab}.  It is true that the
moduli spaces of Higgs bundles and of flat connections carry different
complex structures, but thanks to the correspondence theorem, they are
diffeomorphic, and hence interchangeable topologically.  What is
important for our purposes is that the Higgs moduli space carries a
holomorphic action of the complex torus $\cx$.  This allows the
topology of the moduli space, in some sense, to be determined from
that of the fixed-point set.

Most of the results in this paper are valid for bundles of arbitrary
rank.  We state them this way in the hope that they may hold
independent interest.  We especially have in mind the construction in
\S\ref{limit} of the classifying space of the gauge group as a direct
limit of Higgs bundles.

We originally believed (and claimed in an early version of this paper)
that the entire argument extended to arbitrary rank.  However, a key
lemma, \re{mmm}, is definitely false for rank $>2$.  Nevertheless,
the main theorem, \re{nnn}, turns out to hold true for rank $>2$:
this was established in a paper of Markman \cite{mark2} not long after
this paper was originally written.

{\em Outline of the paper.}  
Section \ref{flat} gives statements of the main results, cast in terms
of flat connections.  Section \ref{review} reviews the correspondence
theorem relating flat connections to Higgs bundles.  The next two
sections give some necessary background on Higgs bundles: section
\ref{p}\ is about their deformation theory, while section
\ref{universe} is about the existence and uniqueness of universal
families.  Section \ref{equiv} then shows how the statements of
\S\ref{flat} will follow from the corresponding (and more general)
statements for Higgs bundles.  From then on Higgs bundles are used
exclusively.

Section \ref{state} re-states the main theorem, on the
generation of the rational cohomology ring by universal classes, in
terms of Higgs bundles.  The proof occupies the remainder of the paper.
It consists of four sections.  Sections \ref{finite} and
\ref{infinite} describe how families of Higgs bundles can be
decomposed into strata where the Harder-Narasimhan filtrations have
fixed type.  Section \ref{finite} is about finite-dimensional
algebraic families and is largely parallel to the work of Shatz
\cite{sha}.  Section \ref{infinite} transfers this stratification to
the setting of an infinite-dimensional space acted on by the gauge
group: a story familiar to readers of Atiyah-Bott \cite{ab}.  

In \S\ref{limit}, however, our approach takes a different turn.  The
moduli spaces under examination are nested in one another, and their
direct limit is shown to have the homotopy type of the classifying
space of the gauge group.  This implies that the cohomology of the
direct limit is indeed generated by universal classes.  It now
suffices to show that this surjects onto the cohomology of the
original moduli spaces.  This is proved in section \ref{surj} by
algebro-geometric methods.  It is here that the rank 2 hypothesis
becomes necessary.

{\em Notation and conventions.}  Throughout the paper, $\rs$ denotes
the compact surface, or smooth projective curve, of genus $g$ over
which we work.  Its cohomology has the usual generators $e_1, \dots,
e_{2g} \in H^1$, and $\si = e_j e_{j+g} \in H^2$.  The Jacobian of
degree $d$ line bundles on $\rs$ is denoted $\Jac^d \rs$; if $d=0$, we
write simply $\Jac \rs$.  Likewise, the $d$th symmetric product of
$\rs$ is denoted $\Sym^d \rs$.  The letters $\H$ and $\M$ denote
moduli spaces over $\rs$, respectively, of $\GL{r,\C}$-connections of
central constant curvature on a bundle of rank $r$ and degree $d$, and
of connections in $\H$ with fixed determinant.  The expressions $\H_n$
and $\M_n$ denote moduli spaces over $\rs$, respectively, of Higgs
bundles of rank $r$ and degree $d$ with values in $K(n) = K \otimes
\co(np)$, and of Higgs bundles $(E, \phi) \in \M_n$ having $\La^r E$
isomorphic to a fixed line bundle $\Xi$ and $\tr \phi = 0$.  Groups
are denoted $T = \cx$, $\Ga = \Sp{2g, \Z}$, and $\Sig = \Z_2^{2g}$.

All cohomology is with rational coefficients unless otherwise stated.
Our methods do not obviously produce any information on the
integral cohomology of the moduli spaces.

We do {\em not} assume $g \geq 2$: the moduli spaces $\M_0$ and $\H_0$
are trivial or empty if $g=0$ or $1$, but $\M_n$ and $\H_n$
for $n>0$ are not so trivial, and they play an important role even for
understanding $g \geq 2$.  

{\em Acknowledgements.} We are grateful to the Mathe\-mat\-isches
Forschungs\-institut Ober\-wolfach for its kind hospitality during
three weeks in the summer of 1998, when the ideas for the present
paper took shape.  We also wish to thank Jeff Brock, Nigel Hitchin,
and Peter Newstead for their advice, help, and encouragement.  The first
author thanks the Institute for Advanced Study and the Miller
Institute for Basic Research in Science for their hospitality and
support during the academic years 1998--99 and 1999--2000,
respectively.  Finally, the second author thanks the Institute for
Theoretical Physics in Santa Barbara, where much of this paper was
written.

\bit{Statement of results in terms of flat connections}
\label{flat}

In one respect the summary given in the introduction is slightly
inaccurate.  The space we study is not exactly the moduli space of
{\em flat} connections on a compact surface.  For that space is
generally singular due to the presence of reducible connections.  This
problem is circumvented by shifting attention to connections of
constant central curvature whose degrees are coprime to their ranks.
The bait-and-switch is perhaps regrettable, but it is standard
practice in the subject.

So let $\rs$ be a surface of genus $g$.  The fundamental group of
$\rs$ has presentation
$$\pi_1(\rs) = 
\frac{\langle a_1, b_1, \dots, a_g, b_g \rangle}
{\prod_{j=1}^g 
a^{\vphantom{-1}}_j b^{\vphantom{-1}}_j a_j^{-1} b_j^{-1}}.
$$

Let $G$ be the non-compact group $\GL{r,\C}$.  A flat $G$-connection
on $\rs$ is determined by a representation $\pi_1 (\rs) \to G$.  So if
$\mu: G^{2g} \to G$ is defined by
$$\mu(A_1, B_1, \dots, A_g, B_g) = \prod_{j=1}^g 
A^{\vphantom{-1}}_j B^{\vphantom{-1}}_j A_j^{-1}
B_j^{-1},$$
then any element of $\mu^{-1}(I)$ defines a flat
$G$-connection.  The quotient $\mu^{-1}(I)/G$, where $G$ acts by
conjugation on all factors, therefore parametrizes all flat
$G$-connections modulo gauge equivalence.  However, as mentioned
above, this is generally a singular space.  An exception is when $r=1$,
but then it is nothing but the complex torus $\Tor = (\cx)^{2g}$.  

Instead of struggling with the singularities, choose any integer $d$
coprime to $r$, and consider the space $\H= \mu^{-1}(e^{2 \pi i d/r}
I) / G$.  This is a non-compact complex manifold --- indeed, a smooth
affine variety \cite{gold,gun} --- and
will be our main object of study.  It parametrizes gauge equivalence
classes of $G$-connections on $\rs$ of constant central curvature
$di\omega I$, where $\omega$ is a 2-form on $\rs$ chosen so that
$\int_\rs \omega = 2 \pi/r$, and $I$ is the $r \times r$ identity
matrix.  Indeed, such a connection is again determined by its
holonomies $(A_j,B_j)$ around $a_j$ and $b_j$, subject to the
constraint that
$$\prod_{j=1}^g 
A^{\vphantom{-1}}_j B^{\vphantom{-1}}_j A_j^{-1} B_j^{-1} 
= \exp \int_\rs di\omega I
= e^{2 \pi di /r} I.$$

Alternatively, $\H$ may be regarded as the space of flat
$G$-connections on $\rs \sans p$ having holonomy $e^{2 \pi di/r} I$
around $p$, modulo gauge.  This has the advantage that no choice of
$\omega$ is necessary, but it is less compatible with the Higgs bundle
interpretation coming up. 

The subject of this paper is the ring structure of the rational
cohomology of $\H$.  To begin this study, some sources of cohomology
classes on $\H$ are needed.

The simplest thing to do is pull back the generators of $H^*(\Tor)$ by
the obvious determinant map $\det: \H \to \Tor$.  This produces
classes $\ep_1, \dots, \ep_{2g} \in H^1(\H)$, but they are not very
interesting.  In fact, the subring they generate can be split off in
the following way.

Let $\M$ be the space constructed in the same way as $\H$, but with
$\SL{r,\C}$ substituted for $\GL{r,\C}$.  Certainly $\M$ is a subspace
of $\H$; but also, scalar multiplication induces a map $\Tor \times \M
\to \H$ which is easily seen to be a free quotient by the abelian
group $\Sig = \Z_r^{2g}$.  According to a theorem of Grothendieck
\cite{toh,mac}, the rational cohomology of such a quotient satisfies
$$H^*(\H) = H^*(\Tor \times \M)^\Sig,$$
where the right-hand side denotes the $\Sig$-invariant part.
Now $\Sig$ acts on $\Tor$ by scalar multiplications, so it
acts trivially on cohomology and hence as rings
$$H^*(\H) = H^*(\Tor) \otimes H^*(\M)^\Sig.$$
Furthermore the composition of the $\Sig$-quotient with the
determinant is the map $\Tor \times \M \to \Tor$ given simply by
projecting to $\Tor/\Sig = \Tor$.  The subring of $H^*(\H)$ generated
by $\ep_1, \dots, \ep_{2g}$ is therefore nothing but the first factor
of the tensor product.

To define more interesting cohomology classes on $\H$, construct a
principal bundle over $\H \times \rs$ as follows.  Let $\overline{G} =
\PGL{r,\C}$.  Any $\rho \in \mu^{-1}(e^{2 \pi di/r} I)$ induces a
well-defined homomorphism $\pi_1(\rs) \to \overline{G}$.  Let
$\tilde{\rs}$ be the universal cover of $\rs$, which is acted on by
$\pi_1(\rs)$ via deck transformations.  There is then a free action of
$\pi_1(\rs) \times G$ on $\overline{G} \times \mu^{-1}(e^{2 \pi di
/r}I) \times \tilde{\rs}$ given by
$$(p,g) \cdot (h,\rho,x) = (\overline{g} \rho(p) h, \overline{g} \rho
\overline{g}^{-1}, p\cdot x),$$
where $\overline{g}$ denotes the image
of $g$ in $\overline{G}$.  The quotient is the desired principal
$\overline{G}$-bundle.  Like any principal $\overline{G}$-bundle, it
has characteristic classes $\bar{c}_2, \dots, \bar{c}_r$, where
$\bar{c}_i \in H^{2i}(\H \times \rs$).  In terms of formal Chern roots
$\xi_k$, $\bar{c}_i$ can be described as the $i$th elementary
symmetric polynomial in the $\xi_k - \zeta$, where $\zeta$ is the
average of the $\xi_k$.  

Now let $\si \in H^2(\rs)$ be the fundamental
cohomology class, and let $e_1, \dots, e_{2g}$ be the basis
of $H^1(\rs)$ Poincar\'e dual to $a_1, \dots, a_g, b_1, \dots, b_g$.  
In terms of these, each of the characteristic classes 
has a K\"unneth decomposition
$$c_i = \al_i \si + \be_i + \sum_{j=1}^{2g} \psi_{i,j} e_j, $$
defining classes $\al_i \in H^{2i-2}(\H)$, $\be_i \in H^{2i}(\H)$, and
$\psi_{i,j} \in H^{2i-1}(\H)$.  The pullback
of these classes to $\Tor \times \M$, by the way, is easily seen to
come entirely from $H^*(\M)^\Sig$.

It is convenient to refer to the entire collection of classes $\al_i$,
$\be_i$, $\psi_{i,j}$, and $\ep_j$ as the {\em universal classes}.

Now specialize, for the remainder of this section, to the case $r=2$.
Then $d$ must be odd, so that $\H = \mu^{-1}(-I)/G$.  The main result
of this paper is then the following.

\bs{Theorem}
\label{nnn}
When $r=2$, the rational cohomology ring of $\H$ 
is generated by the universal classes.  
\es

Equivalently, $\H^*(\M)^\Sig$ is generated by the classes 
$\al = \half \al_2 \in H^2(\M)$, 
$\be = - \ratio{1}{4} \be_2 \in H^4(\M)$, and
$\psi_j = \psi_{2,j} \in H^3(\M)$ for $j=1,\dots,2g$.  
(These normalizations are by now standard in the literature.) 
It is worth mentioning that $H^*(\M)$ is generally
{\em not\/} entirely $\Sig$-invariant \cite[7.6]{h}, hence is not
generated by the universal classes.  However, by the aforementioned
theorem of Grothendieck, $H^*(\M)^\Sig$ is the rational cohomology of
$\M / \Sig$, which is a component of the moduli space of flat
$\overline{G}$-connections on $\rs$.

For completeness, we recount here the main result of our companion
paper \cite{ht} giving all the relations between these generators.
In light of the discussion earlier it suffices to work with
$H^*(\M)^\Sig$.  
Let $\Ga = \Sp{2g,\Z}$.  The action of diffeomorphisms on $\rs$ will
be shown to induce an action of $\Ga$ on $H^*(\M)^\Sig$, fixing $\al$
and $\be$ and acting as the standard representation on the span $V =
H^3(\M)$ of the $\psi_j$.  Thus $\ga = -2 \sum_{j=1}^g \psi_j
\psi_{j+g}$ is a $\Ga$-invariant element of $H^6(\M)$.  Let
$\La^n_0(\psi)$ be the kernel of the natural map $\La^n V \to
\La^{2g+2-n} V$ given by the wedge product with $\ga^{g+1-n}$, or
equivalently, the $\Ga$-invariant complement of $\ga \La^{n-2} V$ in
$\La^n V$.  

For any $g, n \geq 0$, let $I^g_n$ be the ideal within the polynomial
ring $\Q[\al,\be,\ga]$ generated by $\ga^{g+1}$ and the polynomials
$$\rho^c_{r,s,t} 
= \sum_{i=0}^{\min(c,r,s)} \, (c-i)! \, \frac{\al^{r-i}}{(r-i)!} \,
\frac{\be^{s-i}}{(s-i)!} \,\frac{(2 \ga)^{t+i}}{i!},$$
where $c=r+3s+2t-2g+2-n$,
for all $r,s,t \geq 0$ such that $r+3s+3t > 3g-3+n$ and $r+2s+2t
\geq 2g-2+n$.

The main result of our companion paper is then the following.

\bs{Theorem}
\label{dd}
If the rank $r=2$, then as an algebra acted on by $\Ga$, 
$$H^*(\M)^\Sig = \bigoplus_{n=0}^g
\La^n_0(\psi)\otimes\Q[\al, \be, \ga]/I^{g-n}_n.$$
\es

Together, the two main theorems completely describe the ring structure
of $H^*(\H)$ when $r=2$.  They do not completely describe
$H^*(\M)$ because of the classes not invariant under $\Sig$.  However,
these form a relatively minor and simple part of the cohomology, and can
be dealt with by hand; this will be carried out in a forthcoming paper
\cite{thad3}. 

The main theorems will be proved in the language not of flat
connections but rather of Higgs bundles.  Indeed, it proved most
convenient to deduce them from more general results applying to an
infinite sequence of spaces of Higgs bundles, of which $\H$ is only
the first.  We shall next review the definition of Higgs bundles, and
the correspondence theorem relating them to flat connections. 

\bit{Higgs bundles}
\label{review}

A major advance in the study of these representation spaces was made
by Hitchin \cite{h} and Simpson \cite{simp}, who discovered that
they can alternatively be viewed as moduli spaces of holomorphic
objects.  So now, and for the remainder of the paper, let $\rs$ be a
smooth complex projective curve of genus $g$. 

A {\em Higgs bundle\/} or {\em Higgs pair} on $\rs$ with values in a holomorphic line bundle
$L$ is a pair $(E,\phi)$, where $E$ is a holomorphic vector bundle
over $\rs$, and $\phi$, called a {\em Higgs field}, is any element of
$H^0(\End E \otimes L)$.  Its {\em slope} is the rational number $\deg
E / \rk E$.  A holomorphic subbundle $F \subset E$ is {\em
$\phi$-invariant} if $\phi(F) \subset F \otimes L$.  A Higgs bundle
is {\em semistable} if $\slope F \leq \slope E$ for all proper
$\phi$-invariant holomorphic subbundles $F \subset E$, and {\em
stable} if this inequality is always strict.

For example, a pair of the form $(E,0)$ is stable if and only if the
bundle $E$ is stable.

We will be concerned entirely with the case when the line bundle $L$
is $K(n) = K \otimes \co(np)$, where $K$ is the canonical bundle of
$\rs$, $p$ is a distinguished point in $\rs$, and $n \geq 0$.

The moduli space of Higgs bundles with values in $K$ was constructed
by Hitchin \cite{h} and Simpson \cite{simp}, and generalized to an
arbitrary line bundle by Nitsure \cite{nit}.  Their work implies the
following.  

\bs{Theorem}
For fixed rank $r$, degree $d$ coprime to $r$, and $n\geq 0$, there
exists a moduli space $\H_n$ of Higgs bundles with values in $K(n)$,
which is a smooth quasi-projective variety of dimension
$r^2(2g-2+n)+2$.  For a fixed holomorphic line bundle $\Xi$ of degree
$d$, the locus $\M_n$ where $\La^r E \cong \Xi$ and $\tr \phi = 0$ is
a smooth subvariety of dimension $(r^2-1)(2g-2+n)$.  \es

In the case $n=0$, Higgs bundles are related to connections of
constant central curvature in the following way.  Suppose that $\rs$
is equipped with a K\"ahler metric, and let $\omega$ be the K\"ahler
form, again normalized so that $\int_\rs \omega = 2\pi/r$.  Then
Hitchin showed the following.

\bs{Theorem}
\label{v}
Suppose that $r$ and $d$ are coprime and that $n=0$.  Then a Higgs
bundle $(E, \phi)$ is stable if and only if it admits a Hermitian
metric so that the metric connection $A$ satisfies the equation $F_A +
[\phi,\phi^*] = di\omega I$.  This metric is unique up to rescaling
and depends smoothly on $(E, \phi)$. \es

Here $F_A \in \Om^2(\End E)$ is the curvature, the Higgs field $\phi$
is regarded as a section of $\Om^{0,1}(\End E)$, and $I$ is the
identity in $\End E$.  For such a connection $A$, an easy calculation
shows that the $\GL{r,\C}$ connection $A+\phi+\phi^*$ has constant
central curvature $di\omega I$.  
Hence there is a natural smooth map from the space $\H_0$ of Higgs
bundles to the space $\H$ discussed in the previous section.  

In fact, this map is a diffeomorphism, as is the restriction $\M_0 \to
\M$.  The inverse map is provided by a result of Corlette \cite{cor}
and Donaldson \cite{don}.

\bs{Theorem}
\label{zz}
Any $\GL{r,\C}$ connection on $\rs$ with constant central curvature
$di \omega I$ is gauge equivalent to one of the form $A + \phi +
\phi^*$, where $\dbar_A \phi = 0$ and $F_A + [\phi,\phi^*] = di\omega
I$.  \es

Both $\H_0$ and $\H$ carry natural complex structures, but these are
{\em not} identified by the diffeomorphism.  Rather, they are
different members of the family of complex structures which comprises
a {\em hyperk\"ahler structure} on the moduli space.

Our approach does not use this hyperk\"ahler structure.  Indeed, it
will be shown in \re{yy} that the cohomology classes defined above in
terms of flat connections can also be obtained from the universal
family of Higgs bundles.  From then on the flat connection point of
view will vanish, and the moduli space will be regarded exclusively as
a Higgs space. \tag

As pointed out by Simpson \cite{ubi}, the moduli space of Higgs
bundles actually retracts onto a highly singular Lagrangian
subvariety, the {\em nilpotent cone}.  Therefore our results could be
viewed as describing the cohomology ring of the nilpotent cone.
However, this seems to be only a curiosity and is not relevant to our
approach.  

The advantage of the Higgs moduli space is that it admits a
holomorphic action of the group $T = \cx$, given simply by $\la \cdot
(E, \phi) = (E, \la \phi)$. This of course fixes all stable pairs of
the form $(E,0)$, which are parametrized by the moduli space of stable
bundles of rank $r$ and degree $d$.  But the fixed-point set has other
components as well, and they will play a crucial role in what follows.

\bit{Deformation theory of Higgs pairs}
\label{p}

This section and the next summarize, mostly without proof, some basic
facts about Higgs pairs that will be needed later on.  The omitted
proofs are entirely straightforward, along the lines of Markman
\cite[7.3]{mark}, Welters \cite{welt} or the second author
\cite[2.1]{pair}.  In the rank 2 case, some details are worked out in
the first author's Ph.D. thesis \cite{thesis}. 

The deformation space of a holomorphic bundle $E$ is well known to be
$H^1 \End E$; that of a Higgs pair $(E, \phi)$ is similar, but
involves hypercohomology.

Let $(E,\phi)$ be a Higgs pair, and let $\END (E, \phi)$ denote the two-term complex on $\rs$ 
$$\End E \stackrel{[\phantom{\phi},\phi]}{\llrow} \End E \otimes
K(n).$$

\bs{Proposition}
\label{bbb}
The space of infinitesimal deformations of $(E, \phi)$ is the first
hypercohomology group $\Hyp^1 \END (E,\phi)$.  The space of
endomorphisms of $E$ preserving $\phi$ is $\Hyp^0 \END (E,\phi)$.  \fp
\es

Similarly, let $\HOM \lp (E', \phi'), (E, \phi) \rp$ denote the
complex
$$ \Hom (E', E) \lrow \Hom (E', E) \otimes K(n)$$
given by $\psi \mapsto \psi \phi' - \phi \psi$.

\bs{Proposition}
The space of homomorphisms $E' \to E$ intertwining $\phi$ with $\phi'$
is the zeroth hypercohomology group $\Hyp^0 \HOM \lp (E', \phi'), (E,
\phi) \rp$.  The space of extensions of $(E',\phi')$ by $(E, \phi)$ is
$\Hyp^1 \HOM \lp (E', \phi'), (E, \phi) \rp$.  \fp \es

Here an {\em extension} of one Higgs pair by another is a Higgs pair
$(E'', \phi'')$ and a short exact sequence
$$ 0 \lrow E \lrow E'' \lrow E' \lrow 0$$
such that $\phi''$ restricts to $\phi$ on $E$ and projects to $\phi'$
on $E'$.  

One more variation on the theme will be needed in \S\ref{finite}.  Let
$(E, \phi)$ be a Higgs pair containing a flag of $\phi$-invariant
subbundles.  (In practice this flag will always be the {\em
Harder-Narasimhan filtration} defined in \S\ref{finite}.)  Let
$\End' E$ be the subbundle of $\End E$ consisting of endomorphisms
fixing the flag, and let $\END' (E, \phi)$ be the two-term complex
$$\End' E \stackrel{[\phantom{\phi},\phi]}{\llrow} \End' E \otimes
K(n).$$

\bs{Proposition} 
\label{o}
The space of infinitesimal deformations of the Higgs pair $(E,\phi)$
together with the $\phi$-invariant flag is $\Hyp^1 \END' (E, \phi)$. \fp
\es

\bit{Universal families of stable Higgs pairs}
\label{universe}

\bs{Proposition}
\label{m}
Let $(E, \phi)$ and $(E', \phi')$ be stable Higgs pairs with $\slope
E' \geq \slope E$.  Then the dimension of $\Hyp^0 \HOM \lp (E', \phi'), (E,
\phi) \rp$ is $1$ if $(E', \phi')$ and $(E, \phi)$ are
isomorphic, and $0$ otherwise. \fp
\es

In particular, the space of endomorphisms $\Hyp^0 \END (E, \phi)$
consists only of scalar multiplications.

\bs{Proposition}
\label{cc}
Let $(\bee, \bphi)$ and $(\bee', \bphi') \to X \times \rs$ be families
of stable Higgs pairs parametrized by $X$ such that for all $x \in X$,
$(\bee, \bphi)_x \cong (\bee', \bphi')_x$.  Then there exists a line
bundle $L \to X$ such that $(\bee, \bphi) \cong (\bee' \otimes \pi_1^*L,
\bphi)$.  In particular, $\Pj \bee$ and $\End \bee$ are canonical.
\es

\pf.  By  \re{m}, $\Hyp^0 \HOM \lp (\bee', \bphi')_x,
(\bee, \bphi)_x \rp$ is 1-dimensional for all $x$.  Hence the
hyperdirect image $(\R^0 \pi_1)_* \HOM \lp (\bee', \bphi'),
(\bee, \bphi) \rp$ is a line bundle $L \to X$.  It is then easy to
construct the desired isomorphism. \fp

It is clear from the proof that the above proposition holds true not
only for algebraic families of Higgs pairs, but even for smooth
families, that is, $C^\infty$ bundles $(\bee, \bphi) \to X \times \rs$
for any smooth parameter space $X$, endowed with a partial holomorphic
structure in the $\rs$-directions.  

\bs{Lemma} \label{z}
Let $(\bee,\bphi)$ be a family of stable Higgs pairs parametrized by
$X$, and let $\cx$ act on $X$.  If there are two liftings of the
action to $\bee$ so that the induced action on $\bphi$ is $\Ad(\la)
\bphi = \la^{-1}\bphi$, then one is the tensor product of the other
with an action of $\cx$ on a trivial line bundle. In particular, there
are canonical $\cx$-actions on $\Pj \bee$ and $\End \bee$. \es

\pf.  Compose one lifting with the inverse of the other.  This gives a
lifting of the trivial action on $\H_n$ to $\bee$ which preserves
$\bphi$.  By \re{m}, this acts on each fiber via scalar
multiplications.  \fp

\bs{Theorem}
\label{y}
There exists a universal family $(\bee, \bphi)$ over $\H_n \times
\rs$, and a lifting of the $\cx$-action on $\H_n$ to $\bee$ 
whose induced action on $\bphi$ is $\Ad(\la) \bphi = \la^{-1}\bphi$.
\es

That is, $\H_n$ is a fine moduli space for the Higgs bundles of degree
$d$ and rank $r$ with values in $K(n)$.

\pf. This follows in a standard way, cf.\ Newstead \cite{mpos}, from
the geometric invariant theory construction of $\H_n$ due to Nitsure
\cite{nit}.  Alternatively, the universal pair can be constructed
gauge-theoretically just as in Atiyah-Bott \cite[\S9]{ab}. In the rank
2 case, both methods are explained in detail by the first author
\cite[5.3]{hausel2} \cite[5.2.3]{thesis}. \fp

\bit{Equivalence of the two sets of universal classes}
\label{equiv}

The main results were stated in \S\ref{flat} for the moduli space of
flat connections $\H$.  But what is actually proved is similar
results for the Higgs moduli spaces $\H_n$.  To show that these imply
the statements of \S\ref{flat} in the case $n=0$, it suffices to
check that the relevant cohomology classes correspond under the
diffeomorphism $\H_0 \to \H$.  

Let $(\bee, \bphi)$ be a universal family on $\H_n \times \rs$.  There
is a morphism $\H_n\to \Jac^d \rs$ given by $(E, \phi) \mapsto \La^r
E$, so the generators of $H^*(\Jac^d \rs)$ pull back to classes
$\ep_1, \dots, \ep_{2g} \in H^1(\H_n)$.  Also, let $c_2, \dots, c_r$
be the characteristic classes of $\Pj \bee$.  These are elements of
rational cohomology; they can be regarded as the Chern classes of the
tensor product of $\bee$ with a formal $r$th root of $\La^r \bee^*$.

Each of these classes has a K\"unneth decomposition
$$c_i = \al_i \si + \be_i + \sum_{j=1}^{2g} \psi_{i,j} e_j, $$
defining classes $\al_i \in H^{2i-2}(\H)$, $\be_i \in H^{2i}(\H)$, and
$\psi_{i,j} \in H^{2i-1}(\H)$.  The entire collection of classes
$\al_i$, $\be_i$, $\psi_{i,j}$, and $\ep_j$ will be referred to as the
{\em universal classes}.  What requires proof is then the following.

\bs{Theorem}
\label{yy}
When $n=0$, these classes correspond under the diffeomorphism $\H_0
\to \H$ to their counterparts defined in \S\ref{flat}.
\es

\pf.  First, there is a diagram
$$\begin{array}{ccc}
\H & \stackrel{\det}{\lrow} & \Tor \\
\leftdown{\cong} & & \down{} \\
\H_0 & \stackrel{\La^r}{\lrow} & \Jac^d \rs
\end{array}$$
relating the map $\det$ of \S\ref{flat} to the map $(E, \phi) \mapsto
\La^r E$ mentioned above.  It is easy to see from \re{zz} 
that this commutes, and
hence that the classes $\ep_j \in \H$ correspond under the
diffeomorphism.

To show that the higher-degree classes $\al_i$, $\be_i$, and
$\psi_{i,j}$ correspond, it suffices to show that the principal
$\PGL{r,\C}$-bundle associated to $\Pj \bee$ corresponds under the
diffeomorphism 
$\H_0 \to \H$ to the principal bundle 
$$\frac{\overline{G} \times S \times \tilde{\rs}}{\pi_1(\rs) \times
  G}$$
of \S\ref{flat}, where $S = \mu^{-1}(e^{2 \pi id/r}I)$.
(Recall that $G = \GL{r,\C}$, $\overline{G} = \PGL{r, \C}$, and
$\tilde{\rs}$ is the universal cover of $\rs$.)  In fact, we will
construct a principal $G$-bundle $R$ over $\H_0$ and show that the
principal $\overline{G}$-bundle $U$ associated to the pullback of $\Pj
\bee$ to $R \times \rs$ is $G$-equivariantly isomorphic to
$$V= \frac{\overline{G} \times S \times
  \tilde{\rs}}{\pi_1(\rs)}.$$

Let $R$ be the total space of the principal $\overline{G}$-bundle over
$\H_0$ associated to $\Pj \bee |_{\H_0 \times \{ p \} }$.  Then $R$
parametrizes stable pairs $(E, \phi)$ equipped with a frame for the
fiber $E_p$, up to rescaling.  On the other hand, $S$ parametrizes
connections of constant central curvature, together with a frame for
the fiber at a base point, up to rescaling.  The diffeomorphism $\H_0
\to \H$ therefore lifts to a $G$-equivariant diffeomorphism $R \to S$.

Let ${\mathbf F}$ be the pullback to $R \times \rs$ of $\bee$.  Then
$\Pj {\mathbf F}$ admits a natural $G$-action lifting that on $R$, and
it is canonically trivialized on $R \times \{ p \}$.  Moreover, by
\re{v} ${\mathbf F}$ admits a Hermitian metric so that the restriction
of the metric connection ${\mathbf A}$ to each slice $\{ (E, \phi) \}
\times \rs$ satisfies the self-duality equation.  Hence ${\mathbf A} +
\bphi + \bphi^*$ determines a $G$-connection on ${\mathbf F}$ whose
restriction to each slice has constant central curvature.  In
particular, the associated $\overline{G}$-connection on the associated
$\overline{G}$-bundle $U$ is flat on each slice.

On the other hand, the bundle $V$ over $S \times \rs$ defined above is
trivial on $S \times \{ p \}$, and carries a flat connection on each
slice $\{ r \} \times \rs$, which has the same holonomy as the one
just mentioned and is preserved by the action of $G$.  These flat
connections can be used to extend the isomorphism $U|_{R \times \{ p
  \} } \cong V|_{S \times \{ p \} }$ of trivial bundles with
$G$-action to a $G$-equivariant isomorphism $U \cong V$ lifting the
$G$-equivariant diffeomorphism $R \to S$.  \fp

This marks the last appearance of flat connections in our story.  From
now on it is all about Higgs bundles.

\bit{Statement of the generation theorem}
\label{state}

Let $\ep_j$, $\al_i$, $\be_i$, $\psi_{i,j}$ be the universal classes,
defined in \S\ref{equiv}, on the Higgs moduli space $\H_n$.  The goal
of the paper will be to prove this, its main result.

\bs{Theorem}
\label{w}
The rational cohomology ring of $\H_n$ is generated by the universal
classes. 
\es 

The proof of this generation theorem has several parts, with quite
different flavors.  

First, we study the stratification of families of Higgs bundles
according to their Harder-Narasimhan type.  \S\ref{finite} is devoted
to finite-dimensional families, and \S\ref{infinite} to an
infinite-dimen\-sion\-al family analogous to that of Atiyah-Bott
\cite{ab}.  The aim is to show that the strata are smooth of the
expected dimension, but this turns out to be true only in a stable
sense.  We need to consider not only a single $\H_n$, but the chain of
inclusions $\H_n \emb \H_{n+1}$.

We are therefore led to consider in \S\ref{limit} the direct limit
$\H_\infty$ of the $\H_n$, and to show that its cohomology is
generated by universal classes.  Indeed, topological arguments show
that it has the homotopy type of the classifying space of the gauge
group, and the generation then follows from a theorem of Atiyah-Bott
\cite{ab}.

Having done this, we then show in \S\ref{surj} that, when $r=2$, the
cohomology of $\H_\infty$ surjects on that of $\H_n$ for every $n$,
and hence in particular on that of $\H$ itself.  This part of the
proof is algebro-geometric in nature.

\bit{The finite-dimensional stratification}
\label{finite}

We wish to adopt the point of view taken by Atiyah-Bott \cite{ab}, in
which the objects of interest --- for us, Higgs pairs --- are
parametrized by an infinite-dimensional, contractible space.  The
whole space will be divided into strata on which the level of
instability is in some sense constant.  So we will first study the
analogue of this stratification in finite-dimensional algebraic
families, then transfer it to our infinite-dimensional setting.
Except for some subtleties surrounding smoothness, the
results of \S\S\ref{finite} and \ref{infinite} are mostly analogous to
those of Shatz \cite{sha} and Atiyah-Bott \cite{ab}; readers familiar
with those papers may be willing to skip directly to \S\ref{limit}.

Let $(E, \phi)$ be a Higgs pair with values in a line bundle $L$.  A
filtration by $\phi$-invariant subbundles
$$0 = E^0 \subset E^1 \subset \cdots \subset E^\elll = E$$
is said to be a
{\em Harder-Narasimhan filtration} (hereinafter HN filtration) if the
pairs $(F^i,\phi^i)$ are semistable with slope 
strictly decreasing in $i$, where $F^i = E^i /E^{i-1}$ and $\phi^i$ is
induced by $\phi$.

\bs{Proposition}
Any Higgs pair defined over any field of characteristic 0 possesses a
unique HN filtration.
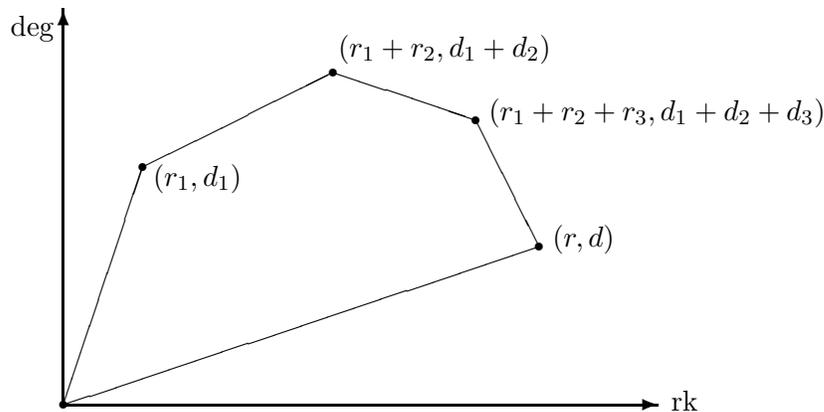
\es

\pf.  The analogous statement for bundles without a Higgs field is
proved by Harder-Narasimhan \cite{hn} and Shatz \cite{sha}.  The
proof for Higgs pairs is entirely parallel: just substitute
$\phi$-invariant subbundles for ordinary subbundles everywhere in
either proof.  Shatz assumes that the ground field is algebraically
closed, but his proof of this theorem does not require it.  \fp

For a given $(E, \phi)$, the {\em type} $\mu$ is the $\elll$-tuple
$(r_1,d_1), \dots, (r_\elll, d_\elll)$ of ranks and degrees of the $F^i$
appearing in its HN filtration.  For example, the type of a semistable
pair is the 1-tuple $(r,d)$.

Since the slope $d_i/r_i$ is strictly decreasing, the pairs $(0,0),
(r_1,d_1), (r_1+r_2,d_1+d_2), \dots, \linebreak[2] (r,d)$ consisting
of partial sums form the vertices of a convex polygon $\Pol(\mu)$ in
$\R^2$, as shown.

\begin{figure}
\begin{center}
\begin{picture}(225,150)
\thicklines
\put(0,0){\vector(0,1){150}}
\put(-20,140){\small deg}
\put(0,0){\vector(1,0){225}}
\put(230,-2){\small rk}
\thinlines
\put(0,0){\circle*{3}}
\put(0,0){\line(1,3){30}}
\put(30,90){\circle*{3}}
\put(33,83){\small $(r_1,d_1)$}
\put(30,90){\line(2,1){72}}
\put(102,126){\circle*{3}}
\put(103,132){\small $(r_1+r_2,d_1+d_2)$}
\put(102,126){\line(3,-1){54}}
\put(156,108){\circle*{3}}
\put(160,108){\small $(r_1+r_2+r_3,d_1+d_2+d_3)$}
\put(156,108){\line(1,-2){24}}
\put(180,60){\circle*{3}}
\put(184,60){\small $(r,d)$}
\put(0,0){\line(3,1){180}}
\end{picture}
\end{center}
\end{figure}

Let $S$ be a scheme of finite type over $\C$, $(\bee, \bphi) \to S
\times \rs$ a family, parametrized by $S$, of Higgs pairs on $\rs$
with values in $L$.

\bs{Proposition} 
\label{h}
The set 
$\{ s \in S \st \mbox{$(\bee, \bphi)_s$ is semistable}\}$
is open in $S$.
\es

\pf.  See Nitsure \cite{nit}.  \fp

For any $\mu$, let $S^\mu$ be the set of those $s \in S$ such that
$(\bee, \bphi)_s$ has type $\mu$.  Recall that a {\em constructible}
subset is a finite union of locally closed sets in the Zariski
topology.

\bs{Lemma}
\label{j}
For any family of Higgs pairs over $S$ and any type $\mu$, $S^\mu$ is
a constructible set.  Moreover, $S^\mu \neq \emptyset$ for only
finitely many $\mu$, and each $S^\mu$ is covered by constructible sets
where the HN filtration varies algebraically, that is,
it determines a filtration of $\bee$ by $\bphi$-invariant subbundles.
\es

\pf.  Without loss of generality assume $S$ is irreducible.  Given a
family $(\bee, \bphi) \to S \times \rs$, let $(\bee, \bphi)_\xi$
be the fiber over the generic point $\xi \in S$.  This is a Higgs pair
defined over the function field $\C(S)$.  It therefore has a HN
filtration, of some type $\mu$, by $\bphi_\xi$-invariant subbundles.
By Lemma 5 of Shatz \cite{sha} there exists an open $U \subset S$ such
that this filtration extends to a filtration of $\bee |_{U \times
\rs}$ by subbundles $\bee^i$.  They are $\bphi$-invariant, since
this is a closed condition and the closure of $\xi$ is all of $S$.

On the other hand, \re{h} implies that, since the quotient pairs
$({\bf F}^i, \bphi^i)_\xi$ are semistable, after restricting to a
smaller $U$ if necessary, $({\bf F}^i, \bphi^i)_s$ are also semistable
for all $s \in U$.  Hence this is the HN filtration at
every $s \in U$, so $U \subset S^\mu$.  

Now pass to $S \sans U$ and proceed by induction on the dimension of
the parameter space.  \fp

Following Shatz \cite{sha}, define a partial ordering on the set of
types by declaring $\mu \leq \nu$ if $\Pol(\mu) \subset \Pol(\nu)$.
Then let $S^{\geq \mu} = \bigcup_{\nu \geq \mu} S^\nu$.

\bs{Theorem} 
\label{i}
For $S$ and $\mu$ as above, $S^{\geq \mu} \subset S$ is closed.
\es

The proof requires the following lemma.

\bs{Lemma}
\label{l}
Let $(E, \phi)$ be a Higgs pair of type $\mu$, and let $F \subset E$
be a $\phi$-invariant subbundle.  Then $(\rk F, \deg F) \in
\Pol(\mu)$.
\es

\pf.  The analogous statement without a Higgs field is Theorem 2 of
Shatz \cite{sha}, and the proof of this is entirely parallel.  One
simply has to note that since the filtration and the subbundle $F$
are $\phi$-invariant, so are the subsheaves $E^i \cap F$ and $E^i
\vee F$ considered by Shatz.  \fp 

\pf\/ of \re{i}.  If $F \subset E$ is any inclusion of torsion-free
sheaves, define $\hat{F} \subset E$ to be the inverse image under the
projection $E \to E/F$ of the torsion subsheaf.  Then $\hat{F}$ and
$E/\hat{F}$ are torsion-free, $F = \hat{F}$ on the locus where $E/F$
is torsion-free, and $F \mapsto \hat{F}$ preserves inclusions of
subsheaves of $E$.

By \re{j} $S^{\geq \mu}$ can be regarded as a
reduced subscheme of $S$.  To show that it is closed, by the valuative
criterion \cite[II 4.7]{har} it suffices to show that, if $X$ is any
smooth curve, and $f: X \to S$ any morphism taking a nonempty open set
into $S^\mu$, then $f(X) \subset S^\mu$.

The generic point of $X$ maps to one of the constructible sets named
in \re{j}, where the HN filtrations are parametrized by subbundles.
Hence there is a non-empty open set $V \subset X$ such that the
restriction of $(f \! \times \! 1)^* \, \bee$ to $V \times \rs$ is
filtered by subbundles restricting over every point in $V$ to the HN
filtration.

Like any coherent subsheaf defined on an open set \cite[II Ex.\ 
5.15(d)]{har}, these bundles extend to coherent subsheaves $\bee^i$ of
$(f \! \times \! 1)^* \, \bee$ over all of $X \times \rs$.  These can
be chosen to remain nested, and as subsheaves of $\bee$ they are of
course torsion-free.  Furthermore, they can be chosen so that
$\bee/\bee^i$ are torsion-free also, by replacing $\bee^i$ with
$\hat{\bee}^i$.

Since torsion-free sheaves on a smooth surface such as $X \times \rs$
are locally free except on a set of codimension 2 \cite[Cor.\ 
2.38]{bob}, it follows that the $\bee^i$ are subbundles except at
finitely many points in the fibers over $X \sans V$.

Now on a smooth curve such as one of these fibers, torsion-free
sheaves are locally free, and so the procedure of the first paragraph
implies the following: every subsheaf of a locally free sheaf is
contained in a subbundle having the same rank and no less degree, with
equality if and only if it was a subbundle to begin with.

When restricted to $\{ x \} \times \rs$ for any $x \in X \sans V$,
then, the nested subsheaves $\bee^i_x$ determine a sequence of
subbundles, which only differ from $\bee^i_x$ at finitely many points
and hence remain nested and $\bphi_x$-invariant, and have the same
rank.   Since the degrees may have risen, they determine a polygon
which contains $\Pol(\mu)$.  By  \re{l}, if the type of $(\bee,
\bphi)_x$ is $\nu_x$, then $\Pol(\nu_x)$ contains this polygon.  Hence
$\nu_x \geq \mu$, so $f(x) \in S^{\geq \mu}$. \fp

As in Atiyah-Bott, the last statement of \re{j} can be refined.

\bs{Proposition}
\label{n}
The HN filtration varies algebraically on all of
$S^\mu$; that is, there exists a filtration of 
$\bee|_{S^\mu \times \rs}$ by subbundles restricting to the
HN filtration on each fiber.
\es

The proof again requires a lemma.

\bs{Lemma}
\label{aaa}
If $E^1$ is the first term in the HN filtration of $(E, \phi)$ and $F
\subset E$ is another $\phi$-invariant subbundle of the same rank and
degree, then $F = E^1$.
\es

\pf.  The corresponding statement for ordinary bundles is a special
case of Lemma 3 of Shatz \cite{sha}.  The proof of this is
again entirely parallel. \fp

\pf\/ of \re{n}.  By  \re{j}, the HN filtrations determine a
constructible subset of the product of Grassmannian bundles $\times_i
\, \mbox{Grass}_{r_i} \bee|_{S^\mu \times \rs}$.  It must be shown
that it is closed.  By the valuative criterion, it suffices to show
that for any morphism $f: X \to S^\mu$, where $X$ is a smooth curve,
the HN filtrations determine a filtration of $(f \! \times \!
1)^* \, \bee$ by subbundles.

As in the proof of \re{i}, a filtration by subbundles does exist over
an open $V \subset X$, and the subbundles extend to nested,
$\bphi$-invariant torsion-free sheaves $\bee^i$ over $X \times \rs$.

The restrictions of these to the fibers over $x \in X \sans V$
generate nested, $\bphi$-invariant subbundles whose ranks and degrees
span a polygon containing $\Pol(\mu)$.  But now since $(\bee,
\bphi)_x$ also has type $\mu$, by  \re{l} this polygon is
contained in $\Pol(\mu)$ as well.  Hence it equals $\Pol(\mu)$, so the
subbundles have degrees equal to those of the subsheaves which
generated them.  They therefore coincide with these subsheaves.

Consequently, the sections of the subsheaf $\bee^i$ span an
$r_i$-dimensional subspace of the fiber of $(f \! \times \! 1)^*
\,\bee$ over every point in $X \times \rs$.  It follows that $\bee^i$
is a subbundle.

Finally, we claim that for any $x \in X \sans V$, the
HN filtration of $\bee_x$ is the restriction of the
$\bee^i$.  First, $\bee^1_x$ is a
$\bphi_x$-invariant subbundle of rank and degree equal to that of the
first term in the HN filtration.  By \re{aaa}
these two subbundles must coincide.  Now pass to
$\bee/\bee^1$ and use induction on the length of the HN filtration to
do the rest. \fp

Recall that $\END'$ refers to the two-term complex, defined in \S\ref{p},
involving endomorphisms fixing a flag of subbundles.  In what
follows, this flag will always be the HN filtration.

\bs{Corollary}
\label{u}
There are deformation maps $T_s S \to \Hyp^1 \END (\bee,
\bphi)_s$ and $T_s S^\mu \to  \Hyp^1 \END' (\bee, \bphi)_s$ so that the
following diagram commutes:
$$\begin{array}{ccc} T_s S^\mu & \lrow & \Hyp^1 \END' (\bee, \bphi)_s \\
\down{} & & \down{} \\
T_s S & \lrow & \Hyp^1 \END (\bee, \bphi)_s. \end{array} $$
\es

\pf.  Follows immediately from \re{bbb}, \re{o} and 
\re{n}.  \fp

\bs{Theorem}
\label{q}
Let $(E, \phi)$ be a Higgs pair with values in $L$.  For $m$ large
enough, $(E, \phi(m))$ belongs to a family, parametrized by a smooth
base $X$, of Higgs pairs with values in $L(m)$ such that the
deformation map $T_{(E, \phi(m))} X \to \Hyp^1 \END (E, \phi(m))$ is an
isomorphism.  \es

\pf.  It suffices to find a smooth $X$ so that the deformation map is
surjective, since then one may restrict to a smooth subvariety
transverse to the kernel.

Choose $k$ large enough that $E(k)$ is generated by its sections and
has $H^1 = 0$.  There is then a surjection $\co^\chi \to E(k)$, where
$\chi$ is the Euler characteristic of $E(k)$.  This represents a point
$q$ in the Quot scheme parametrizing quotients of $\co^\chi$ with
fixed rank and degree.  Let $F$ be the kernel of the map $\co^\chi \to
E(k)$.  The tangent space to the Quot scheme at $q$ is then $H^0 \Hom
(F, E(k))$, and the natural map to the deformation space of $E(k)$ is
the connecting homomorphism in the long exact sequence of
$$0 \lrow \End E \lrow \Hom \lp \co^\chi, E(k)\rp \lrow 
\Hom \lp F, E(k)\rp \lrow \phantom{.}0.$$
This is surjective since $H^1 \Hom \lp \co^\chi, E(k)\rp = \C^\chi
\otimes H^1 (E(k)) = 0$.  For the same reason, $H^1 \Hom \lp F,
E(k)\rp = 0$, so the Quot scheme is smooth at $q$.

Now choose $m$ large enough that $H^1 \lp \End E \otimes L(m) \rp =
0$, and let $X$ be the total space of $\pi_* \lp \End \bee \otimes
L(m) \rp$, where $\bee$ is the tautological quotient on $\Quot \times
\, \rs$, and $\pi$ is projection to $\Quot$.  This is smooth near $q$
since the push-forward is locally free there.  Moreover, there is a
tautological section $\bphi \in H^0 \lp {X \times \rs}, \End \bee
\otimes L(m) \rp$.  That the deformation map is surjective is easily
seen from the diagram
$$
\begin{array}{ccccccc}
H^0 \lp \End E \otimes L(m) \rp & \lrow &
T_{(E, \phi(m))} X & \lrow &
H^0 \Hom \lp F, E(k)\rp & \lrow & 0 \\
\leftdown{=} && \down{} && \down{} && \down{} \\
H^0 \lp \End E \otimes L(m) \rp & \lrow &
\Hyp^1 \END (E, \phi(m)) & \lrow &
H^1 \End E & \lrow & \phantom{.}0. 
\end{array}
$$\fp

Now fix a type $\mu$, and let $(E, \phi)$ be a pair of this type, with
values in $L$.

\bs{Theorem}
\label{t}  
For $m$ large enough, $(E, \phi(m))$ belongs to a family, parametrized
by a smooth base $Y$, of Higgs pairs of type $\mu$ with values in
$L(m)$ such that the deformation map $T_{(E, \phi(m))} Y \to \Hyp^1
\END' (E, \phi(m))$ is an isomorphism.  \es

\pf.  The proof is parallel to that of the previous theorem.  Let $0 =
E^0 \subset \cdots \subset E^\elll = E$ be the HN filtration of $(E,
\phi)$, as usual, and let $E_i = E / E^i$.  Choose $k$ large enough
that every $E_i(k)$ is generated by its sections and every
$H^1(E^i(k)) = 0$.  Then $H^1(E_i(k)) = 0$ as well, and $H^0(E(k)) \to
H^0(E_i(k))$ is surjective.  If $\chi = \dim H^0(E(k))$, a choice of
basis for $H^0(E(k))$ then determines a sequence of quotients
$$ \co^\chi \lrow E(k)
\stackrel{\psi_1}{\lrow} E_1(k) 
\stackrel{\psi_2}{\lrow} E_2(k)
\stackrel{\psi_3}{\lrow} \cdots
\stackrel{\psi_\elll}{\lrow} E_\elll(k) $$
determining a point $q$ in the product of $\elll$ different Quot
schemes.  Let $R$ be the subspace of this product parametrizing flags
of bundles; then $T_q R$ is $\Hyp^0 (C^\cdot)$, where $C^\cdot$ is the
third row in an exact sequence of two-term complexes
$$ \begin{array}{ccc}
0 & & 0 \\
\down{} && \down{} \\
\displaystyle{\medoplus_i} \Hom(E_i, E_i) &\lrow& 
\displaystyle{\medoplus_i} \Hom(E_i,E_{i+1})\\
\down{} && \down{} \\
\displaystyle{\medoplus_i} \Hom(\co^\chi, E_i) &\lrow& 
\displaystyle{\medoplus_i} \Hom(\co^\chi,E_{i+1})\\
\down{} && \down{} \\
\displaystyle{\medoplus_i} \Hom(F_i, E_i) &\lrow& 
\displaystyle{\medoplus_i} \Hom(F_i,E_{i+1}), \\
\down{} && \down{} \\
0 & & 0
\end{array}$$
and $F_i$ is the kernel of the map $\co^\chi \to E_i$.
Now the second row $B^\cdot$ is isomorphic to $\chi$ copies of
$\bigoplus_i E_i \to \bigoplus_i E_{i+1}$, with the map given by
$(b_i) \mapsto (\psi_{i+1} b_i - b_{i+1})$.  Using the long exact
sequence
$$\begin{array}{l}
0 \lrow \Hyp^0 (B^\cdot) \lrow \displaystyle{\medoplus_i} H^0(E_i) \lrow
\displaystyle{\medoplus_i} H^0(E_{i+1}) \\ 
\lrow \Hyp^1 (B^\cdot) \lrow \displaystyle{\medoplus_i}
H^1(E_i) \lrow \displaystyle{\medoplus_i} H^1(E_{i+1})  \\
\lrow \Hyp^2 (B^\cdot) \lrow
0
\end{array} $$
together with $H^1(E_i) = 0$ and the surjectivity of $H^0(E_i) \to
H^0(E_{i+1})$, a descending induction on $i$ shows that
$\Hyp^1(B^\cdot) = \Hyp^2(B^\cdot) = 0$.
Finally, if $A^\cdot$ is the first row, from the short exact sequence
$$0 \lrow \End' E \lrow \medoplus_i \Hom(E_i, E_i) \lrow \medoplus_i
\Hom(E_i,E_{i+1}) \lrow 0$$
we conclude that $\Hyp^* (A^\cdot) = H^*
(\End' E)$ and hence that $\Hyp^0 (C^\cdot) = T_q R$ surjects on
$\Hyp^1(A^\cdot)$ $=$ $H^1(\End' E)$.  Moreover, $\Hyp^2(A^\cdot) =
H^2(\End'E) = 0$.  Hence $\Hyp^1(C^\cdot) = \Hyp^2(C^\cdot) = 0$, so
$R$ is smooth at $q$.

Now choose $m$ large enough that $H^1 \lp \End' E \otimes L(m) \rp =
0$, and let $Y$ be the total space of $\pi_* \lp \End' \bee \otimes
L(m) \rp$, where $\bee$ is the obvious tautological quotient on $R$,
$\End' \bee \subset \End \bee$ is the subbundle of endomorphisms
preserving the flag, and $\pi: R \times \rs \to R$ is the projection.
Then $Y$ is smooth in a neighborhood over $q$, contains a point
representing $(E, \phi(m))$, and has a tautological Higgs field $\bphi
\in H^0 \lp {Y \times \rs}, \End' \bee \otimes L(m) \rp$.  The
deformation map is surjective, as may be seen from the diagram
$$
\begin{array}{ccccccc}
H^0 \lp \End' E \otimes L(m) \rp & \lrow &
T_{(E, \phi(m))} Y & \lrow &
T_q R & \lrow & 0 \\
\leftdown{=} && \down{} && \down{} && \down{} \\
H^0 \lp \End' E \otimes L(m) \rp & \lrow &
\Hyp^1 \END' (E, \phi(m)) & \lrow &
H^1 \End' E & \lrow & \phantom{.}0. 
\end{array}
$$\fp

\bit{The infinite-dimensional stratification}
\label{infinite}

Let $\ce$ be a Hermitian vector bundle over $\rs$ of rank $r$ and
degree $d$.  A rigorous construction of $\H_n$ as an
infinite-dimensional quotient involves connections and sections in
Sobolev spaces associated to $\ce$.  So choose any $k \geq 2$;
Atiyah and Bott prefer $k=2$, but for us as for them, any greater $k$
will also do.  Then, for any Hermitian bundle ${\mathcal V}$ over
$\rs$, denote by $\Om^{p,q}({\mathcal V})$ the Banach space consisting
of sections, of Sobolev class $L^2_{k-p-q}$, of the bundle of
differential forms of type $p,q$ with values in ${\mathcal V}$.  Also
let $\ca$ be the space of holomorphic structures on $\ce$ differing
from a fixed $C^\infty$ one by an element of the Sobolev space
$\Om^{0,1}(\End \ce)$.  We hope the reader will pardon the unorthodox
use of $\Om$ to refer to a Sobolev completion, rather than just the
space of smooth forms.

Define a map 
$$\dbar: \ca \times \Om^{1,0}(\End \ce \otimes K(n)) \to
\Om^{1,1}(\End \ce \otimes K(n))$$
by
$\dbar(E,\phi) = \dbar_E \phi$, and let $\B_n = 
\dbar\vphantom{\partial}^{-1}(0)$.  This
parametrizes all pairs where $\phi$ is holomorphic.

Let $\G$ be the complex gauge group consisting of all complex
automorphisms of $\ce$ of Sobolev class $L^2_k$.  Then $\G$ acts
naturally and smoothly on $\ca$ as shown by Atiyah-Bott, 
and likewise on $\Om^{1,0}(\End \ce \otimes K(n))$ since
$L^2_{k-1}$ is a topological $L^2_k$-module.  The $\G$-action on the
product of these spaces preserves $\B_n$.  

\bs{Proposition}
Every $\G$-orbit in $\B_n$ has a $C^\infty$ representative $(E,\phi)$,
and any two are interchanged by a $C^\infty$ gauge transformation.
The stabilizer of $(E, \phi)$ is the group of holomorphic
automorphisms of $E$ preserving $\phi$.
\es

\pf. According to Lemma 14.8 of Atiyah-Bott, every $\G$-orbit in
$\ca$ contains a $C^\infty$ representative.  If $(E, \phi) \in
\B_n$, so that $\phi$ satisfies the elliptic equation $\dbar_E \phi =
0$, it follows from elliptic regularity that $\phi$ is also
$C^\infty$.  If $(E, \phi)$ and $(E', \phi')$ are pairs in the same
$\G$-orbit, then by Lemma 14.9 of Atiyah-Bott, any gauge
transformation interchanging $E$ and $E'$ is $C^\infty$; hence the
same is true for the pairs.  Finally, if an element of $\G$ preserves
$(E, \phi)$, this means precisely that it preserves $\dbar_E$, hence
is a holomorphic automorphism, and fixes $\phi$. \fp

\bs{Proposition}
\label{r}
Let $(E, \phi) \in \B_n$ be a $C^\infty$ pair.  Then the normal space
to the $\G$-orbit at $(E, \phi)$ is canonically isomorphic to $\Hyp^1
\END (E, \phi)$, and the cokernel of the derivative of $\dbar$ at $(E,
\phi)$ is canonically isomorphic to $\Hyp^2 \END (E, \phi)$.
\es

\pf.  The infinitesimal action of the Lie algebra of $\G$ is the map
$f$, and the derivative of $\dbar$ is the map $g$, in the complex
$$\Om^{0,0}(\End \ce) 
\stackrel{f}{\lrow}
\Om^{0,1}(\End \ce) \oplus \Om^{0,0}(\End \ce \otimes K(n)) 
\stackrel{g}{\lrow}
\Om^{0,1}(\End \ce \otimes K(n)),$$
where $f(a) = (-\dbar_E(a), [a,\phi])$ and $g(b,c) = \dbar_E(c) +
[b,\phi]$.  The symbol sequence of this complex is the
direct sum of those of the Dolbeault complexes of $\End E$ and $\End E
\otimes K(n)$, so it is elliptic.  By elliptic regularity the
cohomology of the complex is then the same as its counterpart where
the Sobolev spaces are replaced by spaces of smooth forms; this is
precisely the Dolbeault hypercohomology of $\END (E, \phi)$.  \fp

\bs{Corollary}
\label{ccc}
For each $(E,\phi) \in \B_n$, $(E,\phi(m))$ is a smooth point of
$\B_{m+n}$ for sufficiently large $m$.
\es

\pf.  Choose $m$ large enough that $H^1(\End E \otimes K(m+n)) =
0$.  Then from the long exact sequence 
$$
\cdots \lrow H^1(\End E) \lrow H^1(\End E \otimes K(m+n))\lrow
\Hyp^2 \END (E, \phi(m)) \lrow 0$$
associated to the hypercohomology of a two-term complex, $\Hyp^2 \END
(E, \phi(m)) = 0$.  Hence $\dbar$ is a submersion at $(E, \phi(m))$,
and the implicit function theorem for Banach manifolds \cite[A3]{dk}
implies that $\B_{m+n} = \dbar\vphantom{\partial}^{-1}(0)$ 
is a smooth embedded Banach
submanifold in a neighborhood of $(E, \phi(m))$.  \fp

Let us now find a slice for the $\G$-action, using the results
of the previous section.  

\bs{Proposition}
\label{s}
Let $(E, \phi)$ be a $C^\infty$ pair in $\B_n$.  
Then for $m$ large enough, there is a $\G$-equivariant submersion
$\G \times U \to \B_{m+n}$ onto a neighborhood of $(E, \phi(m))$, where
$U$ is an open neighborhood of $x$ in the algebraic family of \re{q}.  
\es 

\pf.  \re{q} provides a family $(\bee, \bphi)$ of pairs over
some smooth $X \ni x$ such that $(\bee, \bphi)_x = (E, \phi(m))$ and
the natural map $T_x X \to \Hyp^1 \END (E, \phi(m))$ is an isomorphism.
Choose a Hermitian metric on $\bee$ extending the given one on $E$.
This determines a smooth map $X \to \B_{m+n}$, which by 
\re{r} is transverse to the $\G$-orbit.  It extends to a
$\G$-equivariant map $\G \times X \to \B_{m+n}$ whose derivative is a
surjection at $\G \times \{ x \}$, and hence in some neighborhood $\G
\times U$, thanks to the diagram
$$\begin{array}{ccccccc}
T_e\G &\lrow& T_{(e,x)}(\G \times X) &\lrow& T_x X &\lrow& 0\\
\leftdown{=} && \down{} && \down{} && \down{}  \\
\Om^{0,0}(\End E) &\lrow& T\B_n &\lrow& \Hyp^1 \END (E,\phi) &\lrow& 
\phantom{.}0.
\end{array} $$ \fp

Let $\B_n^\mu$ denote the union of all $\G$-orbits in $\B_n$ whose
$C^\infty$ representatives have type $\mu$ in the sense of
\S\ref{finite}.  In particular, let $\B_n^s$ denote the stable orbits.

\bs{Corollary}
For any $\mu$, $\bigcup_{\nu \geq \mu}
\B_n^\nu$ is closed in $\B_n$.
\es

\pf.  Follows immediately from \re{s} and the
corresponding fact for $X$, \re{i}.  \fp

Let $\overline{\G}$ be the quotient of $\G$ by the central subgroup
$\cx$.  Since by \re{m} the latter is the stabilizer of
all stable pairs, each stable orbit is isomorphic to $\overline{\G}$.

\bs{Corollary}
\label{ddd}
The natural map $\B_n^s \to \H_n$ is a principal
$\overline{\G}$-bundle.
\es

\pf.  The submersions of \re{s} descend to maps $\overline{\G} \times U
\to \B_n$, whose derivatives are isomorphisms on $\overline{\G} \times
\{ x \}$, and which can be made injective by shrinking $U$ if
necessary.  By the inverse function theorem for Banach manifolds
\cite[A1]{dk} these are $\overline{\G}$-equivariant diffeomorphisms
onto their images.  They therefore constitute an atlas of local
trivializations.  \fp

Given a $C^\infty$ pair $(E,\phi)\in \B_n$, define $\End'' E$ by the short exact sequence
$$0 \lrow \End' E \lrow \End E \lrow \End'' E \lrow 0,$$
where $\End' E$, as before, is the subsheaf of $\End E$ preserving the
HN filtration of $(E, \phi)$.
Also let $\END'' (E, \phi)$ be the two-term complex on $\rs$
defined analogously to $\END (E,\phi)$ and $\END' (E,\phi)$.
There is then a short exact sequence of two-term complexes
$$ {\bf 0} \lrow \END' (E,\phi) \lrow \END (E,\phi) 
\lrow \END'' (E,\phi) \lrow {\bf 0}.$$

\bs{Theorem}
\label{ggg}
For any $(E, \phi) \in \B_n$ of type $\mu$, and for $m$ large enough,
$\B_{m+n}^\mu$ is an embedded submanifold of $\B_{m+n}$ near $(E,
\phi(m))$ with normal space canonically isomorphic to $\Hyp^1 \END''
(E, \phi(m))$.  \es

\pf.  By acting with an element of $\G$ if necessary we may assume
that $(E, \phi)$ is $C^\infty$.  

For $m$ large, there was constructed in \re{t} a family
of pairs $(\bee, \bphi)$ of type $\mu$ over a smooth base $Y \ni y$,
having $(\bee, \bphi)_y = (E, \phi(m))$ and $T_y Y \to 
\Hyp^1 \END' (E, \phi(m))$ an isomorphism.  Choose a metric on $\bee$
extending the given one on $E$.  This determines a smooth map $Y \to
\B_{m+n}$, which by \re{r} is transverse to the
$\G$-orbit. 

On an open neighborhood $V$ of $y$ in $Y$, choose a lifting of this
map to the domain of the submersion $\G \times U \to \B_{m+n}$ of
\re{s}.  Projecting this lifting to $U$ gives a map $V \to
U$ whose image consists of pairs of type $\mu$, and hence is contained
in $X^\mu$.  Its derivative at $y$ is the natural map from $T_y Y =
\Hyp^1 \END' (E, \phi(m))$ to $T_x X = \Hyp^1 \END (E, \phi(m))$.

This derivative is injective.  Indeed, the kernel is the image of $\Hyp^0
\END'' (E, \phi(m))$.  But if $0 = E^0 \subset \cdots \subset E^\elll =
E$ is the HN filtration of $(E, \phi)$ as usual, then there is a short
exact sequence 
$$0 \lrow \End'' E^{\elll-1} \lrow \End'' E \lrow \Hom (E^{\elll-1},
E/E^{\elll-1}) \lrow 0$$
and hence a short exact sequence of the
corresponding two-term complexes.  The first hypercohomology $\Hyp^0$
of both of the outer complexes vanishes by an induction on $\elll$
using \re{m}, so $\Hyp^0 \END'' (E, \phi(m)) = 0$ too.

Hence, in a neighborhood of $x$, $X^\mu$ contains an embedded
submanifold with tangent space $\Hyp^1 \END' (E, \phi(m)) \subset
\Hyp^1 \END (E, \phi(m))$.  But by \re{u}, $X^\mu$ is a
subscheme of $X$ with Zariski tangent space contained in $\Hyp^1 \END'
(E, \phi(m))$.  It therefore must be smooth near $x$.  

The inverse image of $\B_{m+n}^\mu$ under the submersion $\G \times U
\to \B_{m+n}$ is $\G \times (U \cap X^\mu)$, so this immediately
implies that $\B_{m+n}^\mu$ is a smoothly embedded submanifold in a
neighborhood of the orbit of $(E, \phi(m))$.  Its normal space at $(E,
\phi(m))$ is the quotient of $\Hyp^1 \END (E, \phi(m))$ by $\Hyp^1
\END' (E, \phi(m))$; by choosing $m$ large enough we may arrange as in
\re{ccc} that
$\Hyp^2 \END' (E, \phi(m)) = 0$, so that this quotient is nothing but
$\Hyp^1 \END'' (E, \phi(m))$. \fp

\bit{The direct limit of Higgs spaces}
\label{limit}

The inclusions $\B_n \emb \B_{n+1}$ make the set of all $\B_n$ into a
directed set.  Let $\B_\infty$ be the direct limit.  It may be
regarded as a set of pairs $(E,\phi)$ as before, but where $\phi$ may
now have a pole of arbitrary finite order at $p$.  Note that for each
type $\mu$, the direct limit of $\B_n^\mu$ is a subset $\B_\infty^\mu$
of $\B_\infty$.  Note also that $\G$ acts naturally on $\B_\infty$ and
that $\B^s_\infty /\G$ is $\H_\infty$, the direct limit of the
$\H_n$.  In another context, $\H_n$ has appeared in the work of
Donagi-Markman \cite{dm}. 

Our aim in this section is to show that $\B^s_\infty$ is contractible.
Essentially, the reason is that $\B_\infty$ is contractible, and the
complement of the stable set has infinite codimension.

Recall that a subspace of a topological space is a {\em deformation
neighborhood retract} (hereinafter DNR) if it is the image of a map
defined on some open neighborhood of itself and homotopic to the
identity.  It is {\em equivariant} if the homotopy is
equivariant for the action of some group.

\bs{Lemma}
\label{eee}
As a subspace of $\H_\infty$, each $\H_n$ is a DNR, and the open sets
which retract can be chosen to be nested.
\es

\pf.  Since $\H_n \emb \H_{n+1}$ is an embedding of finite-dimensional
manifolds, we may choose a tubular neighborhood $U_n^1$ and a
projection $U_n^1 \to \H_n$.  This tubular neighborhood in turn has a
tubular neighborhood $U_n^2$ in $\H_{n+2}$, namely its inverse image
in $U_{n+1}^1$, and so on.  The direct limits $U_n^\infty$ of these
tubular neighborhoods are nested open subsets of $\H_\infty$.  Each is
homeomorphic to a vector bundle over $\H_n$ and hence deformation
retracts onto it.  Indeed, the deformation retraction preserves each
$U_n^i$. \fp

\bs{Theorem}
\label{fff}
As a subspace of $\B_\infty$, each $\B_n^s$ is a $\G$-equivariant DNR,
and the $\G$-invariant open sets which retract can be chosen to be
nested. 
\es

\pf. Let $\pi: \B_\infty^s \to \H_\infty$ be the quotient map, and let
$U_n^i$ be the tubular neighborhoods of the previous proof.  By the
definition of direct limit, it suffices to construct a sequence of
$\G$-equivariant deformation retractions of $\pi^{-1}(U_n^i)$ onto
$\B_n^s = \pi^{-1}(\H_n)$ compatible with the inclusions
$\pi^{-1}(U_n^i) \subset \pi^{-1}(U_n^{i+1})$.  Obviously we would
like to lift the deformation retracts of the previous proof to the
principal $\G$-bundle.

These liftings are guaranteed to exist by the first covering homotopy
theorem.  This asserts that if $F: Y \times [0,1] \to Z$ is a homotopy
with any reasonable domain (including manifolds, but unfortunately not
their direct limits), and if $E$ is a fiber bundle over $Z$, then
$F^*E$ is isomorphic as a fiber bundle to $F^*E|_{Y \times 0} \times
[0,1]$.  Actually, a slight refinement of this result is needed,
namely that the isomorphism can be chosen so as to extend a given one
over a closed DNR $X \subset Y$.  We then apply this refined result to
the case where $X = U_n^i$, $Y=U_n^{i+1}$, and the homotopy is the
retraction of $U_n^i$ on $\H_n$ described above.  It is easy to
construct the desired deformation retraction of $\pi^{-1}(U_n^i)$ from
the resulting isomorphism.

The slight refinement of the first covering homotopy theorem can be
proved by the same argument as the theorem itself, given by Steenrod
\cite[\S11.3]{steen}.  Just choose the atlas for $F^*E$ so that its
restriction to $X \times [0,1]$ is pulled back from an atlas on $X$
using the given isomorphism.  The existence of such an atlas follows
easily from the fact that $X$ is a closed DNR and the ordinary version
of the theorem.  \fp

\bs{Corollary}
\label{jjj}
The quotient map $\pi: \B_\infty^s \to \H_\infty$ is a principal
$\overline{\G}$-bundle.
\es

\pf.  Immediate from \re{ddd} and \re{fff}.  \fp

\bs{Theorem}
\label{hhh}
For all $k \geq 0$, $\pi_k(\B_\infty^s) = 1$.
\es

\pf.  The open subsets of $\B_\infty$ provided by \re{fff}, which
retract onto $\B_n^s$, form a nested open cover of
$\B_\infty^s$.  By compactness, any map $S^k \to \B_\infty^s$, and
any homotopy of such maps, has image contained in one such
neighborhood.  Hence $\pi_k(\B_\infty^s) = \lim_{n \to \infty} \pi_k
(\B_n^s)$.  

So let $f: S^k \to \B_n^s$ be any map.  Now $\B_n$ is contractible
just by retracting it first on $\ca \times \{ 0 \}$, so $f$ certainly
extends to a continuous map $f: D^{k+1} \to \B_n$.  Our task is to
perturb this so that it misses the unstable locus.

For each $x \in D^{k+1}$, by \re{ccc} and \re{ggg} there is some
integer $m_x \geq n$ such that for all $m \geq m_x$, $f(x)$ is a
smooth point of $\B_m$, and the stratum $\B^\mu_m$ containing $f(x)$
is an embedded submanifold at $f(x)$.  Passing to a finite subcover
and taking $m = \max m_x$, we find that $f(D^{k+1})$ maps entirely
into the smooth locus of $\B_m$, and that near its image each stratum
$\B_m^\mu$ is an embedded submanifold of finite codimension.

Yet another application of compactness shows that $f(D^{k+1})$
intersects only a finite number of strata $\B_m^\mu$.  By increasing
$m$ again if necessary we may assume by \re{ggg} that each of these
strata has codimension $> k+1$.  Then, starting with the stratum of
highest codimension and working our way up, we may perturb $f$ so that
it no longer touches that stratum, but so that its value on $S^k$
remains unchanged.  After these perturbations, $f$ will have image
entirely within $\B_m^s$.

Thus for any $f: S^k \to \B_n^s$, the homotopy class of $f$ is killed
by the inclusion in $\B_m^s$ for some $m \geq n$, so $\lim_{n \to
\infty} \pi_k (\B_n^s) = 1$. \fp

For an alternate proof in the rank 2 case, see the first author's
thesis \cite[7.5.1]{thesis}.

\bs{Theorem}
\label{iii}
The space $\B_\infty^s$ is contractible.
\es

\pf.  A theorem of Whitehead \cite[10.28]{sw} asserts that if $X$ is a
CW-space --- that is, a space homotopy equivalent to a CW-complex ---
whose homotopy groups all vanish, then it is contractible.  In light
of \re{hhh} above, it therefore suffices to show that
$\B_\infty^s$ is a CW-space.  

Consider the fiber bundle
\beq \label{g}
\B_\infty^s \lrow \frac{\B_\infty^s \times E
  \G}{\G} \lrow B\G.
\eeq
We will show that this is a fibration whose total space and base space
are CW-spaces.  It then follows from Corollary 13 of Stasheff \cite{sta}
that the fiber is a CW-space.  

Note that $\G$ acts on $\B_\infty^s$ with stabilizer $\cx$ and
quotient $\H_\infty$.  There is therefore a fiber bundle
$$B\cx \lrow \frac{\B_\infty^s \times E \G}{\G} \lrow \H_\infty,$$
which is just the associated bundle to the principal
$\overline{\G}$-bundle $\B_\infty^s \to \H_\infty$.  The base of this
fiber bundle, being a direct limit of manifolds, is metrizable and
hence paracompact by Stone's theorem \cite[6-4.3]{munk}; hence the
fiber bundle is a fibration \cite[I 7.13]{white}.  Moreover, the base
is a CW-space, as is the fiber.  Proposition 0 of Stasheff \cite{sta}
then implies that the total space is a CW-space.

According to Proposition 2.4 of Atiyah-Bott \cite{ab}, $B\G$ is a
component of the space of maps from $\rs$ to $BG$, with the
compact-open topology.  Since the domain is a compact metric space and
the range is a CW-complex, by Corollary 2 of Milnor \cite{mil2} the
space of maps is a CW-space.  Moreover, since $BG$ is metrizable and
$\rs$ is compact, the space of maps $B\G$ is metrizable, hence
paracompact.  The fiber bundle \re{g} is therefore a fibration.  \fp

\bs{Corollary}
\label{kkk}
The space $\H_\infty$ is homotopy equivalent to $B\overline{\G}$.
\es

\pf.  By \re{jjj} and \re{iii}, there is a principal
$\overline{\G}$-bundle on $\H_\infty$ with contractible total space.  \fp

Those who dislike the appearance of infinite-dimensional,
gauge-theoretic methods in the last two sections may wish to reflect
that it is no doubt possible to replace every infinite-dimensional
construction by a finite-dimensional, algebraic approximation, in the
style of Bifet et al.\ \cite{bgl} or Kirwan \cite{k2}.  In any case,
algebraic geometry will reappear on the scene shortly.

\bit{Surjectivity of the restriction on cohomology}
\label{surj}

Let $\ce$ be as in \S\ref{infinite}.  The pull-back of
$\ce \to \rs$ to the product $\B^s_\infty \times \rs$ is acted on by
$\G$, so the projective bundle $\Pj \ce$ is acted on by
$\overline{\G}$.  It therefore descends to a $\Pj^r$-bundle $\Pj \bee$
over $\H_\infty \times \rs$, whose characteristic classes can be
decomposed as usual into K\"unneth components:
$$\bar{c}_i = \al_i \si + \sum_{j=1}^{2g} \psi_{i,j} e_j + \be_i.$$

Likewise, the natural determinant maps $\H_n \to \Jac^d \rs$ given by
$(E, \phi) \mapsto \La^r E$ extend to $\H_\infty \to \Jac^d \rs$.  Let
$\ep_1, \dots, \ep_{2g}$ be the pull-backs of the standard generators
of $H^1(\Jac^d \rs)$.

It is straightforward to check that the universal classes $\al_i$,
$\be_i$, $\psi_{i,j}$, $\ep_j$ thus defined restrict to their
counterparts on $\H_n$ for $n \geq 0$.

\bs{Theorem}
The rational cohomology ring $H^*(\H_\infty)$ is generated by
these universal classes.
\es

\pf.  According to Atiyah-Bott, $H^*(B\G)$ is generated by universal
classes, which means the following.  First of all, $B\G$ can be
identified with the component of the space of maps $\rs \to B\U{r}$
such that the pull-back of the universal bundle over $B\U{r}$ is
isomorphic to $\ce$.  Atiyah-Bott call this component $\Map_\ce
(\rs,B\U{r})$.  Then, the pull-back of the universal bundle by the
canonical map $\Map_\ce (\rs, B\U{r}) \times \rs \to B\U{r}$ is a
bundle whose Chern classes can be decomposed into K\"unneth components
as usual.  Atiyah-Bott prove \cite[2.20]{ab} that these generate the
ring $H^*(B\G)$.

On the other hand, by \re{kkk} $B\overline{\G}$ can also be identified
with $\H_\infty$.  Hence $B\G$ is a bundle over $\H_\infty$ with fiber
$B\cx = \C \Pj^\infty$.  As explained by Atiyah-Bott, the rational
cohomology of this bundle splits: $H^*(B \G) = H^*(B \overline{\G})
[h]$.  By restricting to a single $\C \Pj^\infty$ fiber, it can be
checked that $\be_1 = r h$ modulo elements of $H^2(B \overline{\G})$;
it may therefore be discarded since we seek only generators of $H^*(B
\overline{\G})$.  Also $\al_1 \in H^0(B\G)$ may be discarded since by
\re{hhh} $B\G$ is connected.  Finally it can be checked that
$\psi_{1,j} = \ep_j$, and that for $i > 1$, the classes $\al_i$,
$\be_i$ and $\psi_{i,j}$ of Atiyah-Bott agree with those defined
above.  (Strictly speaking, they may differ by some lower order terms,
since the characteristic classes of a projective bundle are evaluated
by formally twisting so that $c_1 = 0$.) \fp

Now by \re{eee}, $H_*(\H_\infty)$ is the direct limit of
$H_*(\H_n)$, and hence $H^*(\H_\infty)$ is the inverse limit of
$H^*(\H_n)$.  Consequently, the surjectivity of the restriction map
$H^*(\H_\infty) \to H^*(\H_n)$ for all $k$, and hence the generation
theorem \re{w}, is implied by the following
result, whose proof occupies the remainder of this section.

\bs{Theorem}
\label{lll}
When $r=2$, the restriction $H^*(\H_{n+1}) \to H^*(\H_n)$ is surjective.
\es

Let $\T = \cx$ act on each $\H_n$ by $\la \cdot (E,\phi) = (E, \la
\phi)$.  This action is compatible with the inclusion $\H_n \emb
\H_{n+1}$.  Furthermore, since this is an algebraic action on a smooth
quasi-projective variety, the $\U{1}$-part of the action is
Hamiltonian, and the $\R^\times$ part of the action is the Morse flow
of the moment map.

\bs{Lemma}
For any $(E, \phi) \in \H_n$, there exists a limit $\lim_{\la \to 0}
(E, \la \phi) \in \H_n$.
\es

\pf.  We may regard this as a limit of the downward Morse flow in
$\H_n$.  Note that it need not be simply $(E, 0)$ as this may be
unstable.  Nevertheless, a stable limit always exists; this may be
seen in two ways.

First, one can regard $\H_n$ as a space of solutions $(A, \phi)$ to
the self-duality equations, as Hitchin \cite{h} regards $\H$; the
moment map is then $(A, \phi) \mapsto \| \phi \|^2$, and what we need
to know is that this is proper and bounded below.  The boundedness is
obvious, and the properness is proved following Hitchin's argument for
$\H$ \cite[7.1(i)]{h}.

Alternatively and more algebraically, one can observe, as does Simpson
\cite[Prop.\ 3]{ubi}, that the Hitchin map defined by Nitsure
\cite[6.1]{nit}, taking $\H_n$ holomorphically to a vector space, is
proper and intertwines the $\T$-action on $\H$ with a linear action on
the vector space having positive weights.  A limit must therefore
exist in the zero fiber of the Hitchin map.  \fp

Now any $\U{1}$ moment map whose Morse flows have lower limits is a
{\em perfect\/} Bott-Morse function: see for example Kirwan
\cite[9.1]{kir}.  This means that its Morse inequalities are
equalities.  More explicitly, it means the following.  Let $y_0,
\dots, y_k$ be the critical values of the moment map $\mu: X \to \R$,
and $F_i$ the corresponding critical submanifolds.  Choose real
numbers $x_i$ so that $x_0 < y_0 < x_1 < y_1 < \cdots < y_k < x_{k+1}$.
If $X_i = \mu^{-1}(x_0, x_i)$, then, as for any Bott-Morse function,
there is a homotopy equivalence of pairs $(X_{i+1},X_i) \simeq
(D_i,S_i)$, where $D_i$ is the disc bundle, and $S_i$ the sphere
bundle, associated to the {\em negative normal bundle} of $F_i$, that
is, the bundle of downward Morse flows, cf.\ Milnor \cite{mil}.  For
the function to be {\em perfect\/} means that moreover the connecting
homomorphism vanishes in each long exact sequence
$$ \cdots \lrow H^*(X_{i+1},X_i) \lrow H^*(X_{i+1}) \lrow H^*(X_i)
\lrow \cdots, $$
breaking it up into short exact sequences.

Suppose now that $X$ contains a $\T$-invariant submanifold $Y$ on
which the moment map is again perfect.  Then by induction on $i$,
$H^*(X)$ surjects on $H^*(Y)$ if and only if $H^*(X_{i+1},X_i)$
surjects on $H^*(Y_{i+1},Y_i)$ for all $i$.

We find ourselves in this situation, with $X = \H_{n+1}$ and $Y
=\H_n$.  To prove \re{lll}, it therefore suffices to show that the
relative cohomology of the disc bundle for the downward flow from each
critical submanifold in $\H_{n+1}$ surjects on that of its
intersection with $\H_n$.  This will be true in the case $r=2$.
Indeed, the Thom isomorphism identifies this relative cohomology with
the ordinary cohomology of the critical submanifold itself.  We will
prove, first, that this identification is compatible with the
restriction to $\H_n$, and second, that the latter restriction is
surjective.  Both will follow from the description below of the critical
set (cf.\ Hitchin \cite[7.1]{h}).

\bs{Proposition}
\label{ooo}
The critical submanifolds of the moment map on $\H_n$ are a disjoint union
$$\F_n = \bigsqcup_{j=0}^{g+ \left[ \frac{n-1}{2} \right]} F_n^j,$$
where: 
\begin{enumerate}
\item for $j=0$, the absolute minimum $F_n^0$ of the moment map is the
moduli space of stable bundles of rank $2$ and degree $d$,
parametrizing Higgs bundles $(E,\phi)$ with $\phi = 0$;
\item for $j>0$, $F_n^j = \Jac^{\frac{d+1}{2} - j} \rs \times
\Sym^{2g+n-1-2j} \rs$,
parametrizing Higgs bundles $(E,\phi)$ with $E = L \oplus M$, $\deg L
= \frac{d+1}{2} - j$, and 
$$\phi = \left( \begin{array}{cc} 0&0\\s&0 \end{array} \right) ,$$
where $s \in H^0(KML^{-1}(n))$ vanishes on an effective divisor of
degree $2g+n-1-2j$.
\end{enumerate}
\es

\pf.  The critical points for the moment map of the action of $\U{1}
\subset T$ are exactly the fixed points of $T$.  For any $(E,\phi) \in
\H_n$ fixed by $T$, by \re{y} $T$ acts by automorphisms on the
universal bundle restricted to $\{ (E, \phi) \} \times \rs$, which is
nothing but $E$ itself, and $\la \in T$ takes $\phi$ to $\la \phi$.
If the weights of the action are distinct, this
splits $E$ as a sum of line bundles $L \oplus M$, and $\phi$ is forced
to be of the stated form.  If the weights are not distinct, then $T$
acts by scalars, so $\phi$ is invariant and hence must be $0$.  \fp

\bs{Lemma}
\label{mmm}
The downward flow from $F_n^j = F_{n+1}^j \cap \H_n$ in $\H_{n+1}$ is
wholly contained in $\H_n$.  
\es

\pf.  The statement is vacuous for $j=0$, since at the absolute
minimum there is no downward flow.  Consider then a point $(E, \phi)
\in F_n^j$ for as described in \re{ooo}(b).  According to \re{bbb},
the tangent space to $\H_{n+1}$ at $(E,\phi)$ is the hypercohomology
$\Hyp^1 \END (E,\phi)$, where $\END (E,\phi)$ is the two-term complex 
$$\End E \stackrel{[\phantom{\phi},\phi]}{\llrow} \End E \otimes
K(n+1).$$
Since $E = L \oplus M$ and $\phi$ is strictly lower-triangular, this
breaks up as a direct sum of complexes, which are the weight spaces
for the $T$-action.  The downward flow corresponds to $\Hyp^1$ of the
complex 
$$\Hom(L,M) \lrow 0,$$
which is of course just $H^1(ML^{-1})$.  This is independent of $n$,
and so the downward flow in $\H_{n+1}$ is wholly contained in $\H_n$,
as desired.  \fp

Consequently, we may use the Thom isomorphisms to identify the map of
relative cohomology on the disc bundles with the ordinary restriction
map $H^*(\F_{n+1}) \to H^*(\F_n)$.  To prove \re{lll},
then, it remains to prove the following statement.

\bs{Theorem}
\label{d}
The restriction $H^*(\F_{n+1}) \to H^*(\F_n)$ is surjective.
\es

\pf.  Again, this is vacuous for $j=0$, since $F_{n+1}^0 = F_n^0$.  It
is not much harder for $j>0$, for $F_n^j$ is isomorphic to a product of a
Jacobian and a symmetric product, and the embedding $F_n^j \emb F_{n+1}^j$
corresponds to the identity on the first factor and a map of effective
divisors of the form $D \mapsto D + p$ on the second.  The latter map
is easily seen (for example, from the description of Macdonald
\cite{mac}) to induce a surjection on cohomology.  \fp
  
This completes the proof of \re{lll}, and hence of the generation
theorem \re{w}.

Note, by the way, that the theorem
is {\em false} for the moduli space $\M_0$ consisting of
pairs with fixed determinant line bundle $\La^n E$ and trace-free
$\phi$.  One can see this already from Hitchin's description
\cite[7.6]{h} of its cohomology: the universal classes are
all invariant under the natural action of $\Sig = \Z_2^{2g}$, but
there also exist classes which are not $\Sig$-invariant.  (See also
our companion paper \cite[\S4]{ht}.) \tag

In the cases of ranks 2 and 3, the first author has found an
alternative proof of the generation theorem, which will appear in a
forthcoming paper \cite{mum}.  It uses the vector
bundles over $\H_n$ whose Chern classes above their ranks furnish the
so-called ``Mumford relations''; these are well-known to have the
appropriate dimension in the $\GL{r}$ case, but not in the $\SL{r}$
case.  

The generation theorem does, however, tell us the following about the
cohomology ring of the fixed-determinant moduli space.

\bs{Corollary}
Let $\Xi$ be a fixed line bundle of odd degree, and let
$\M_n$ be the moduli space of Higgs bundles $(E, \phi) \in \H_n$ with
$\La^2 E \cong \Xi$ and $\tr \phi = 0$.  Then $\Sig \cong \Z_2^{2g}$ acts
naturally on $\M_n$, and 
$$H^*(\H_n) = H^*(\Jac \rs) \otimes H^*(\M_n)^\Sig,$$
where $H^*(\Jac \rs)$ is generated by the $\ep_1, \dots, \ep_{2g}$,
and $H^*(\M_n)^\Sig$ by the remaining universal classes.
\es

\pf.  In fact $\H_n$ is the quotient of $T^* \Jac \rs \times \M_n$ by
the free action of $\Sig$; the quotient map is $(L, \psi) \times (E,
\phi) \mapsto (L \otimes E, \psi \id + \phi)$.  A theorem of
Grothendieck referred to earlier \cite{toh} asserts that the rational
cohomology of a quotient by a finite group is the invariant part of
the rational cohomology.  Hence \beqas
H^*(\H_n) & = & H^*(T^* \Jac \rs \times \M_n)^\Sig \\
& = & H^*(\Jac \rs) \otimes H^*(\M_n)^\Sig, \eeqas since $\Sig$ acts
on $T^* \Jac \rs$ only by translations.  The determinant map $\H_n \to
\Jac^d \rs$ lifts to the map $T^* \Jac \rs \times \M_n \to \Jac^d \rs$
given by projection on the first factor followed by an isogeny of
order $r$; hence $\ep_1,\dots, \ep_{2g}$ generate the rational
cohomology of the first factor.  A universal pair on $\H_n$ pulls back
to the tensor product of the Poincar\'e line bundle on $\Jac \rs
\times \rs$ with a universal pair on $\M_n \times \rs$; hence its
projectivization, and consequently the remaining universal classes,
are pulled back from the second factor.  \fp

\end{document}